\documentclass{article}  

\usepackage{algorithm,algorithmic}
\usepackage{amsmath,amssymb,amsthm}
\usepackage{latexsym}
\usepackage{graphicx}
\usepackage{epsfig}
\usepackage[left=35mm, right=35mm, 
            top=35mm, bottom=35mm]{geometry}

\newtheorem{theo}{Theorem}
\newtheorem{prop}[theo]{Proposition}
\newtheorem{lem}[theo]{Lemma}
\newtheorem{defi}[theo]{Definition}
\newtheorem{rem}[theo]{Remark}

\newcommand{\set}[1]{\left\{#1\right\}}
\newcommand{\card}{\mathrm{card}}

\newcommand{\ie}{{\em i.e.}\ }
\newcommand{\eps}{\varepsilon}

\newcommand{\Flip}[3]{\phi(\mathbf{#1},#2,#3)}

\newcommand{\SsChemins}[1]{\mathcal{A}_{#1}}
\newcommand{\Meandres}[1]{\mathcal{M}_{#1}}
\newcommand{\Walls}[1]{\mathcal{W}_{#1}}
\newcommand{\Culmis}[1]{\mathcal{C}_{#1}}
\newcommand{\Exc}[1]{\mathcal{E}_{#1}}

\begin{document}

\title{Random sampling of lattice paths with constraints, via transportation}

\author{Lucas Gerin}

\maketitle

\begin{abstract}
We investigate Monte Carlo Markov Chain (MCMC) procedures for the random 
sampling of some one-dimensional lattice paths with constraints, for various constraints.
We will see that an approach inspired by \emph{optimal transport} 
allows us to efficiently bound the mixing time of the associated Markov chain. The algorithm is robust and easy to
implement, and samples an "almost" uniform path of length $n$ in $n^{3+\eps}$ steps.
This bound makes use of a certain \emph{contraction property} of the Markov chain, and is also used
to derive a bound for the running time of Propp-Wilson's \emph{Coupling From The Past} algorithm.
\end{abstract}



\section{Lattice Paths with Constraints}

Lattice paths arise in several areas in probability and combinatorics, either in their own interest
(as realizations of random walks, or because of their interesting combinatorial properties:
see \cite{Ban} for the latter) or because of fruitful bijections with 
various families of trees, tilings, words.
The problem we discuss here is to efficiently sample uniform (or \emph{almost} uniform)
paths in a family of paths with constraints.

There are several reasons for which one may want to generate uniform samples of lattice paths: 
to make and try conjectures on the behaviour of a large "typical" path, test 
algorithms running on paths (or words, trees,...). 
In view of random sampling, it is often very efficient to make use of the combinatorial structure of the
family of paths under study. In some cases, this yields linear-time (in the length of the path) 
\emph{ad-hoc} algorithms \cite{MBM,Duc}. However, the nature of the constraints makes sometimes impossible
such an approach, and there is a need for robust algorithms that work in lack of 
combinatorial knowledge.  

Luby,Randall and Sinclair \cite{LRS} design a Markov chain that generate
sets of non-intersecting lattice paths. This was motivated by a classical (and simple, see illustrations in \cite{Des,Wilson}) correspondence between dimer configurations on an hexagon, rhombae tilings of this hexagon and families of non-intersecting lattice paths.
As the first step for the analysis of this chain, Wilson \cite{Wilson} introduces a peak/valley
Markov chain (see details below) over some simple lattice paths and obtain sharp bounds for its mixing
time.
We present in this paper a variant of this Markov chain,
which is valid for various constraints and whose analysis is simple. 
It generates an "almost" uniform path of length $n$ in $n^{3+\eps}$ steps, this bound makes use of a certain
\emph{contraction property} of the chain.

Appart from the algorithmic aspect, the peak/valley process seems to have a physical relevancy as a simplified model for the evolution of \emph{quasicrystals} (see a discussion on a related process in the introduction of \cite{Des}). In particular, the mixing time of this Markov seems to have some importance.

\subsection*{Notations}
\begin{figure}[h]
\begin{center}
\includegraphics[width=40mm]{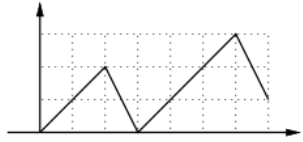}
\caption{The lattice path $S=(1,2,0,1,2,3,1)$ associated with the word $(1,1,-2,1,1,1,-2)$.}
\end{center}
\end{figure}
We fix three integers $n,a,b>0$, and consider the paths of length $n$, with steps $+a/-b$,
that is, the words of $n$ letters taken in the alphabet $\set{a,-b}$. 
Such a word $s=(s_1,s_2,\dots,s_n)$ is identified to the path
$S=(S_1,\dots,S_n):=(s_1,s_1+s_2,\dots,s_1+s_2+\dots +s_n)$. 

\begin{figure}[h]
\begin{center}
\includegraphics[width=65mm]{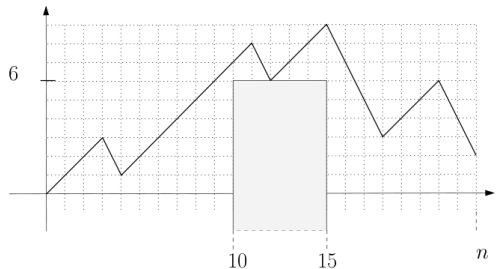}
\caption{A path of steps $+1/-2$, with a wall of height $h=6$ between $i=10$ and $j=15$.}
\label{Fig:CheminMur}
\end{center}
\end{figure}

To illustrate the methods and the results, we focus on some particular sub-families 
$\SsChemins{n}\subset \set{a,-b}^n$:
\begin{enumerate}
\item Discrete \emph{meanders}, denoted by $\Meandres{n}$, which are simply the non-negative paths: $S\in\Meandres{n}$ if for
any $i\leq n$ we have $S_i\geq 0$. This example is mainly illustrative because the combinatorial properties of meanders
make it possible to perform exact sampling very efficiently (an algorithm running in $\mathcal{O}(n^{1+\eps})$
steps is given in \cite{MBM}, an order that we cannot get in the present paper).
\item Paths with \emph{walls}. A path with a wall of height
$h$ between $r$ and $s$ is a path such that
$S_i\geq h$ for any $r\leq i\leq s$ (see Fig. \ref{Fig:CheminMur} for an example).
These are denoted by $\Walls{n}=\Walls{n}(h,r,s)$, they
appear in statistical mechanics as toy models for the analysis of random interfaces
and polymers (see examples in \cite{Walls}).
\item \emph{Excursions}, denoted by $\Exc{n}$, which are non-negative paths such that $S_n=0$. 
In the case $a=b=1$, these correspond to well-parenthesed words and are usually called Dyck words.
In the general case, Duchon \cite{Duc} proposes a rejection algorithm which generates excursions
in linear time.
\item \emph{Culminating paths} of size $n$, denoted further by $\Culmis{n}$, which are non-negative paths whose maximum is attained at the last step: for any $i$ we have $0\leq S_i\leq S_n$. They have been introduced in \cite{MBM}, 
motivated in particular by the analysis of some algorithms in bioinformatics.
\end{enumerate}

\section{Sampling with Markov chains}\label{Sec:Sampling}
We will consider Markov chains in a family $\SsChemins{n}$, where all the probability transitions are symmetric. For a modern introduction to Markov chains, we refer to \cite{Hagg}.
Hence we are given a transition matrix $(p_{i,j})$ of size $|\SsChemins{n}|\times|\SsChemins{n}|$ with 
\begin{align*}
p_{i,j} &=p_{j,i} \mbox{ whenever }i\neq j,\\
p_{i,i} &= 1-\sum_{j\neq i}p_{i,j}.
\end{align*}
\begin{lem}\label{Lem:Unif}
If such a Markov chain is irreducible, then it admits as unique stationary distribution
the uniform distribution $\pi=\pi(\SsChemins{n})$ on $\SsChemins{n}$.
\end{lem}

\begin{proof}
The equality $\pi(i) p_{i,j}= \pi(j) p_{j,i}$ holds for any two vertices $i,j$. This shows that
the probability distribution $\pi$ is reversible for $(p_{i,j})$, and hence stationary. It is unique 
if the chain is irreducible.
\end{proof}

This lemma already provides us with a scheme for sampling an almost uniform path in $\SsChemins{n}$, without knowing much about $\SsChemins{n}$.
To do so, we define a ``flip'' operator on paths, this is an operator
$$
\begin{array}{r c c c}
\phi: & \SsChemins{n}\times \set{1,\dots,n}\times \set{\downarrow,\uparrow}\times \set{+,-} &\to &\SsChemins{n}\\
      & (\mathbf{S},i,\eps,\delta)                                                           &\mapsto & \phi(\mathbf{S},i,\eps,\delta).
\end{array}
$$
When $i\in\set{1,2,\dots,n-1}$ the path $\phi(\mathbf{S},i,\uparrow,\delta)$ is defined as follows : 
if $(s_i,s_{i+1})=(-b,a)=$ \includegraphics[width=7mm]{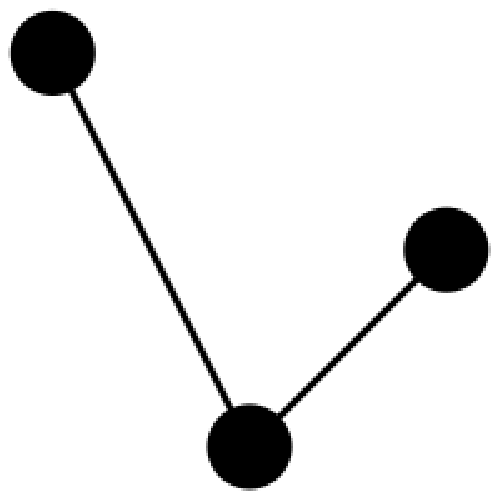}
then these two steps are changed into $(a,-b)=$ \includegraphics[width=7mm]{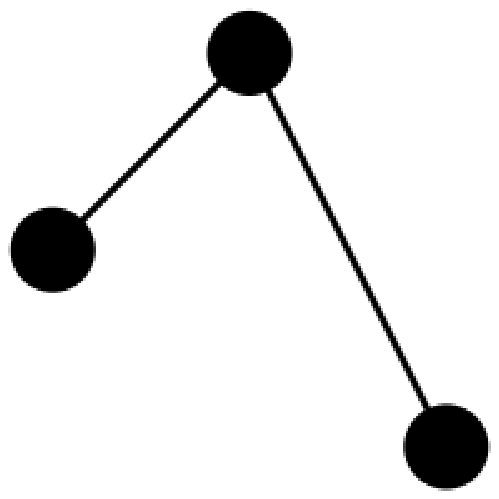}. The $n-2$ other steps remain unchanged. 
If $(s_i,s_{i+1})\neq (-b,a)$ then $\Flip{S}{i}{\uparrow}{\delta}=\mathbf{S}$. Note that in the case 
$i\in\set{1,2,\dots,n-1}$ the value of $\phi$ does not depend on $\delta$.

For the case $i=n$, if $\delta=+$, we define $\Flip{S}{n}{\eps}{\delta}$ as before as if there would be a $+a$ as the 
end if the path. For instance, in the case where $S_n=-b$, the path $\Flip{S}{n}{\uparrow}{+}$, the $n$-th step is turned into $a$.

The path $\Flip{S}{i}{\downarrow}{\delta}$ is defined equally: if $i<n$ and $(s_i,s_{i+1})=$ \includegraphics[width=7mm]{updown.eps},
it turns into \includegraphics[width=7mm]{downup.eps}. When $\delta=-$, one flips as if there would be a
$-b$ at the end of the path.

For culminating paths, we have to take another definition
of $\Flip{S}{n}{\uparrow}{\delta},\Flip{S}{n}{\downarrow}{\delta}$, see Section \ref{Sec:Analysis}.

We are also given a probability distribution $\mathbf{p}=(p_i)_{1\leq i\leq n}$, and we assume that $p_i>0$ for each $i$. 
We will consider a particular sequence $\mathbf{p}$ later on, but at this point we can take
the uniform distribution in $\set{1,\dots,n}$. We describe the algorithm below in Algorithm \ref{Algo:CM}.
\begin{algorithm}
\caption{Approximate sampling of a path in $\SsChemins{n}$}
\label{Algo:CM}
\begin{algorithmic}
\STATE initialize $\mathbf{S}=(+a,+a,+a,\dots,+a)$
\STATE $I_{1},I_{2},\dots\leftarrow$ i.i.d. r.v. with law $\mathbf{p}$
\STATE $\eps_{1},\eps_{2},\dots\leftarrow$ i.i.d. uniform r.v. in $\set{\uparrow,\downarrow}$
\STATE $\delta_{1},\delta_{2},\dots\leftarrow$ i.i.d. uniform r.v. in $\set{+,-}$
\FOR{$t=1$ to $T$}
  \IF{$\Flip{S}{I_t}{\eps_t}{\delta_t}$ is in $\SsChemins{n}$} \STATE $\mathbf{S}\leftarrow  
\Flip{S}{I_t}{\eps_t}{\delta_t}$ 
  \ENDIF
\ENDFOR
\end{algorithmic}
\end{algorithm}

In words, this algorithm performs the Markov chain in $\SsChemins{n}$ with transition matrix
$P=\left(P_{\mathbf{R},\mathbf{S}}\right)_{\mathbf{R},\mathbf{S}\in\SsChemins{n}}$ defined as follows:
$$
\begin{cases}
P_{\mathbf{R},\mathbf{S}}&=p_i/2, \mbox{ if }  \mathbf{S}\neq\mathbf{R}\text{ and } \mathbf{S}=\Flip{R}{i}{\eps}{\delta} \mbox{ for some }i,\eps, \delta\\
P_{\mathbf{R},\mathbf{S}}&=0\text{ otherwise,}\\
P_{\mathbf{R},\mathbf{R}}&=1-\sum_{\mathbf{S}\neq \mathbf{R}} P_{\mathbf{R},\mathbf{S}}.\\
\end{cases}
$$

\begin{prop}
Denote by $S(t)$ the random path obtained after the $t$-th run of the loop in 
Algorithm \ref{Algo:CM}. When $t\to\infty$, the 
sequence $S(t)$ converges in law to the uniform distribution in $\SsChemins{n}$.
Moreover, the execution of Algorithm \ref{Algo:CM} until time $T$ is linear in $T$.
\end{prop}
\begin{proof}
For the first claim, we have to check that the chain is aperiodic and irreducible.
Aperiodicity comes from the (many) loops. Irreducibility will follow from 
Lemma \ref{Lem:Geodesique}.
For the second claim, notice that the time needed for the test "$\Flip{S}{I_t}{\eps_t}$ is in $\SsChemins{n}$" can be considered as constant, since for the families $\Meandres{n}$ and $\Exc{n}$ we only have to 
compare $0,S_i$ while for the family $\Walls{n}$ we only have to compare $S_i$ with the height of the wall
at $i$. For the case of the culminating paths, see below in Section \ref{Sec:Analysis}.
\end{proof}

We now choose the distribution $(p_i)$. Instead of $p_i=1/n$, we will use the probability distribution defined 
by
\begin{equation}\label{Eq:Poids}
p_i:=i(2n-i)\kappa_0  +a\quad (\mbox{ for }i=1,\dots,n),
\end{equation}
where
\begin{align*}
\kappa_0&=\frac{3}{2n^2(n+1)}\\
a       &=1/4n^3.
\end{align*}
We let the reader check that $(p_i)_{i\leq n}$ 
is indeed a probability distribution. 
The reason for which we use this particular distribution
will appear in the proof of Proposition \ref{Lem:Courbure}.
We will then need the following relation: for 
each $1\leq i\leq n-1$,
\begin{equation}\label{Eq:kappa}
p_i-p_{i-1}/2-p_{i+1}/2 = \kappa_0.
\end{equation}

For Algorithm \ref{Algo:CM} to be efficient, we need to know how $S(T)$ is close in law to
$\pi$.
This question is related to the spectral properties of the matrix $P$. In particular, the speed of convergence is governed by the spectral gap (\ie $1-\lambda$, where $\lambda$ is the largest of the modulus of the eigenvalues different from one,
see \cite{Mix} for example), but this quantity is not known in general.
Some geometrical methods \cite{Dia} allow to bound from below $1-\lambda$, but they assume a precise knowledge of the 
structure of the graph defined by the chain $P$. It seems that such results do not apply here.

Instead, we will study the metric properties of the chain $P$ with respect to a natural distance on 
$\SsChemins{n}$, and show that
it satisfies a certain \emph{contraction property}. 

\subsection{The variant of Algorithm \ref{Algo:CM} for culminating paths}\label{Sec:Analysis}

Unchanged, our Markov chain $P$ cannot generate culminating paths since the path $\mathbf{S}=(a,a,\dots,a)$
would then be isolated: it has no peak/valley 
and $\Flip{S}{n}{\downarrow}{-}=(a,a,\dots,-b)$ which is not culminating. 

Thus we propose a slight modification for the family $\Culmis{n}$. We only change the definition
of $\Flip{S}{i}{\eps}{\delta}$ when $i=n$ (it won't depend on $\delta$).
Since the maximum is reached at $n$, the  $\lceil b/a\rceil +1$ last steps are necessarily
$$
(a,a,\dots,a) \mbox{ or } (-b,a,\dots,a).
$$
We thus define $\Flip{S}{n}{\uparrow}{\delta}$ as the path obtained by changing the $\lceil b/a\rceil +1$ last steps into $(a,a,\dots,a)$ (regardless of their initial values in $\mathbf{S}$)
and $\Flip{S}{n}{\downarrow}{\delta}$ as the path obtained by changing the $\lceil b/a\rceil +1$ last steps into $(-b,a,\dots,a)$.

Notice that despite this change the execution time of each loop of Algorithm \ref{Algo:CM} is still a $\mathcal{O}(1)$:
\begin{itemize}
\item If $I_t< n$, the time needed for the test "$\Flip{S}{I_t}{\eps_t}{\delta_t}$ is in $\SsChemins{n}$" can be considered as constant, since we only have to compare $0,S_i,S_n$.
\item If $I_t=n$, the new value $S_n$ is compared with the maximum of $S$, which can be done in $\mathcal{O}(n)$. Fortunately, this occurs with probability $p_n=\mathcal{O}(1/n)$, so that the time-complexity of each loop is, on average, a $\mathcal{O}(1)$.
\end{itemize}

\section{Error estimates with contraction}\label{Sec:Ricci}

Going back to a more general setting, we consider a Markov chain in a finite set $V$, endowed with a metric $d$.
For a vertice $x\in V$ and a transition matrix $P$, we denote by $P\delta_x$ (resp. $P^t\delta_x$) 
the law of the Markov chain associated with $P$ at time $1$ (resp. $t$), when starting from $x$.
For $x,y\in V$, the main assumption made on $P$ is that there is a coupling 
between $P\delta_x,P\delta_y$ (that is, a random variable
$(X_1,Y_1)$ with $X_1\stackrel{\mbox{law}}=P\delta_x,Y_1\stackrel{\mbox{law}}=P\delta_y$) such that 
\begin{equation}\label{Eq:Courbure}
\mathbb{E}\left[d(X_1,Y_1)\right]\leq (1-\kappa)d(x,y),
\end{equation}
for some $\kappa >0$, which is called the \emph{Ricci curvature} of the chain, by analogy with the Ricci curvature in 
differential geometry\footnote{The Ricci curvature is actually the largest positive number such 
that \eqref{Eq:Courbure} holds, for all the couplings of $P\delta_x,P\delta_y$ ; here we should rather say 
that Ricci curvature is larger than $\kappa$.}. 
If the inequality holds, then it implies (\cite{Mix},p.189) that $P$ admits a unique stationary 
measure $\pi$ 
and that, for any $x$,
\begin{equation}\label{Eq:Mixing}
\parallel P^t\delta_x -\pi\parallel_{\mathrm{TV}} \leq (1-\kappa)^t\mathrm{diam}(V),
\end{equation}
where $\mathrm{diam}(V)$ is the diameter of the graph with vertices $V$ induced by the Markov chain.
The notation $\parallel .\parallel_{\mathrm{TV}}$ stands, as usual, for the \emph{Total Variation} distance
over the probability distributions on $V$ defined by
$$
\parallel \mu_1 -\mu_2\parallel_{\mathrm{TV}}:= \sup_{A\subset V} \left|\mu_1(A)-\mu_2(A)\right|.
$$ 
Hence, a positive Ricci curvature gives the exponential convergence to the stationary measure, with an exact (\ie 
\eqref{Eq:Mixing} is non-asymptotic in $t$) bound. In many situations, a smart choice for the coupling between $X_1,X_2$ gives a sharp rate of convergence in eq. \eqref{Eq:Mixing} (see some striking
examples in \cite{Olli}).

\subsection{Metric properties of $P$}
To apply the Ricci curvature machinery, we endow each $\SsChemins{n}$ with the $L^1$-distance
$$
d_1(S,S')=\frac{1}{a+b}\sum_{i=0}^n |S_i-S_i'|.
$$
(Notice that $|S_i-S_i'|$ is always a multiple of $a+b$.) For our purpose, it is fundamental that this metric
space is \emph{geodesic}.
\begin{defi}
A Markov chain $P$ in a finite set $V$ is said to be \emph{geodesic} with respect to the distance $d$ 
on $V$ if for any two $x,y\in V$ with $d(x,y)=k$, there exist $k+1$ vertices $x_0=x,x_1,\dots,x_k=y$ in $V$ such that for each $i$
\begin{itemize}
\item $d(x_i,x_{i+1})=1$ ;
\item $x_i$ and $x_{i+1}$ are neighbours in the Markov chain $P$ (\ie $P(x_i,x_{i+1})>0$).
\end{itemize}
This implies in particular that $P$ is irreducible and that the diameter of $P$ is smaller 
than $\max_{x,y}d(x,y)$.
\end{defi}
\begin{lem}\label{Lem:Geodesique}
For each family $\Culmis{n}$,$\Walls{n}$,$\Exc{n}$,$\Meandres{n}$ the Markov chain of Algorithm \ref{Algo:CM} is geodesic with respect to $d_1$. 
\end{lem}
\begin{proof}[Proof of Lemma \ref{Lem:Geodesique}]
The proof goes by induction on $k$. We fix $S\neq T$ (and denote by $s,t$ the corresponding words) ; we 
want to decrease $d_1(S,T)$ by one, by applying the operator $\phi$ with proper $i,\eps$. 
We denote by $i_0\in\set{1,\dots,n}$ the first index for which $S\neq T$.
For instance we have $T_{i_0}=S_{i_0}+a+b$. Let $j$ be the position of the left-most peak in $T$ in $\set{i_0+1,i_0+2,\dots,n}$, if such a peak exists.
Then $S':=\Flip{T}{j}{\downarrow}{\delta}$ is also in $\SsChemins{n}$: it is immediate for the 
families $\Meandres{n},\Walls{n},\Culmis{n},\Exc{n}$.
We have $d_1(S,S')=k-1$.

If there is no peak in $T$ after $i_0$, then $(t_{i_0+1},t_{i_0+2},\dots,t_n)=(a,a,\dots,a)$. Hence we try to 
increase the final steps of $S$ by one. 
To do so, we choose $S':=\Flip{S}{n}{\uparrow}{\delta}$ if $S\neq \Flip{S}{n}{\uparrow}{\delta}$, 
or $S'=\Flip{S}{j}{\uparrow}{\delta}$  where $j$ is the position of the right-most $-b$ otherwise
(we choose the right-most one to ensure that $\Flip{S}{j}{\uparrow}{\delta}$ remains culminating in the case where
$\SsChemins{n}=\Culmis{n}$.).

\end{proof}

For meanders, excursions and walls, we will show that the Ricci curvature of $P$ with respect to
the distance $d_1$ is (at least) of order $1/n^3$. 

\begin{prop}\label{Lem:Courbure}
For the three families $\Meandres{n},\Exc{n},\Walls{n}$, the Ricci curvature of the associated Markov chain,
with weights $(p_i)$ defined as in \eqref{Eq:Poids}, is larger than $\kappa_0$.
\end{prop}
\begin{proof}[Proof of Proposition \ref{Lem:Courbure}]
Fix $\mathbf{S},\mathbf{T}$ in $\SsChemins{n}\in\set{\Meandres{n},\Exc{n},\Walls{n}}$, 
we first assume that $\mathbf{S},\mathbf{T}$ are neighbours, for instance $\mathbf{T}=\Flip{S}{i}{\uparrow}$
for some $i$.
\begin{center}
\includegraphics[width=35mm]{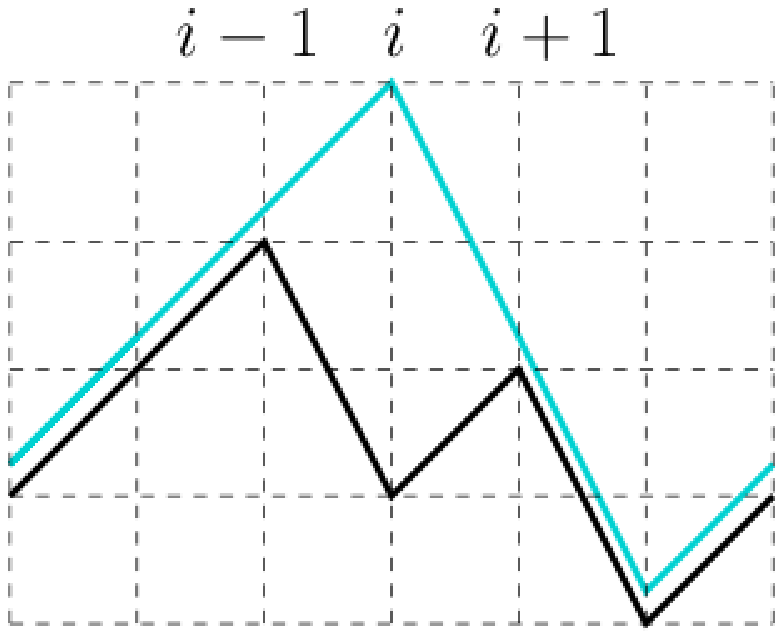}
\end{center}
Let $(\mathbf{S}^1,\mathbf{S}^2)$ be the random variable in $\SsChemins{n}\times\SsChemins{n}$
whose law is defined by
$$
(\mathbf{S}^1,\mathbf{S}^2)\stackrel{\mbox{(law)}}{=} \left(\phi(\mathbf{S},\mathcal{I},{E}),\phi(\mathbf{T},\mathcal{I},E)\right),
$$
where $\mathcal{I}$ is a r.v. taking values in $\set{1,\dots,n}$ with distribution $\mathbf{p}$ and
$E$ is uniform in $\set{\uparrow,\downarrow}$.
In other words, we run one loop of Algorithm \ref{Algo:CM} simultaneously on both paths.

We want to show that $\mathbf{S}^1,\mathbf{S}^2$ are, on average, closer than $\mathbf{S},\mathbf{T}$.
Different cases may occur, depending on $\mathcal{I}$ and on the index $i$ where $\mathbf{S},\mathbf{T}$ differ.

\noindent{\bf \underline{Case 1. $i=1,2,\dots,n-2$.}\ \ } 
\par
\noindent{\bf \underline{Case 1a. $\mathcal{I}=i$.}\ \ } 
This occurs with probability $p_i$ and, no matter the value of $E$, we have $\mathbf{S}^1=\mathbf{S}^2$.

\par
\noindent{\bf \underline{Case 1b. $\mathcal{I}=i-1$ or $i+1$.}\ \ } 
We consider the case $i-1$. Since $\mathbf{S}$ and $\mathbf{T}$ coincide 
everywhere but in $i$, we necessarily have one of these two cases:
\begin{itemize}
\item there is a peak in $\mathbf{S}$ at $i-1$ and neither a peak nor a valley in $\mathbf{T}$ at $i-1$ (as in the figure on the right) ;
\item there is a valley in $\mathbf{T}$ at $i-1$ and neither a peak nor a valley in $\mathbf{S}$ at $i-1$.
\end{itemize}
In the first case for instance, then we may have $d_1(\mathbf{S}^1,\mathbf{S}^2)=2$ if $E=\downarrow$, while
the distance remains unchanged if $E=\uparrow$. The case $\mathcal{I}=i+1$ is identical.
This shows that with a probability smaller than 
$p_{i-1}/2+p_{i+1}/2$ we have $d_1(\mathbf{S}^1,\mathbf{S}^2)=2$.

\par
\noindent{\bf \underline{Case 1c. $\mathcal{I}\neq i-1,i,i+1$ and $\mathcal{I}\neq n$.}\ \ } 
In this case, $\mathbf{S}$ and $\mathbf{T}$ are possibly modified in $\mathcal{I}$, but if there is
a modification it occurs in both paths. It is immediate since for the families
$\Meandres{n}$,$\Walls{n}$ and $\Exc{n}$ since the constraints are local.

\par
\noindent{\bf \underline{Case 2. $i=n-1$.}\ \ } 
In this case, it is easy to check that, because of our definition of $\Flip{S}{n}{\eps}{\delta}$, 
we have
$$
\mathbb{E}\left[d_1(\mathbf{S}^1,\mathbf{S}^2)\right]
\leq 1-p_{n-1}+p_{n-2}/2+p_{n}/2 =1-\kappa_0.
$$

\par
\noindent{\bf \underline{Case 3. $i=n$.}\ \ } 
We have
$$
\mathbb{E}\left[d_1(\mathbf{S}^1,\mathbf{S}^2)\right]
\leq 1+p_{n-1}/2-p_{n}/2=1-\kappa_0.
$$

\vspace{3mm}

Thus, we have proven that when $\mathbf{S},\mathbf{T}$ only differ at $i$
\begin{align}
\mathbb{E}\left[d_1(\mathbf{S}^1,\mathbf{S}^2)\right]
&\leq 2\times(p_{i-1}/2+p_{i+1}/2) +0\times p_i+1\times(1-p_i-p_{i-1}/2-p_{i+1}/2)\label{Eq:Accroissement}\\
&\leq (1-\kappa_0)\times 1=(1-\kappa_0)d_1(S,T)\notag.
\end{align}
What makes Ricci curvature very useful is that if this inequality 
holds for pairs of neighbours then it holds for 
any pair, as noticed in \cite{Bub}. Indeed, take $k+1$ paths $S_0=S,S_1,\dots,S_k=T$ as in Lemma 
\ref{Lem:Geodesique} and apply the triangular inequality for $d_1$:
\begin{align*}
\mathbb{E}\left[d_1(\phi(S,\mathcal{I},{E}),\phi(T,\mathcal{I},E))\right]
&\leq \sum_{i=0}^{k-1} \mathbb{E}\left[d_1(\phi(S_i,\mathcal{I},{E}),\phi(S_{i+1},\mathcal{I},E))\right]\\
&\leq (1-\kappa_0)k=(1-\kappa_0)d_1(S,T).
\end{align*}
\end{proof}
\begin{rem}
It is easy to exhibit some $S,T$ such that ineq. \eqref{Eq:Accroissement} is in fact
an equality. In 
the case where $p_i=1/n$, this equality reads $\mathbb{E}\left[d_1(\mathbf{S}^1,\mathbf{S}^2)\right]=d_1(S,T)$,
and we cannot obtain a positive Ricci curvature (though this does not
prove that there is not another coupling or another distance for which we could 
get a $\kappa >0$ in the case $p_i=1/n$.).
\end{rem}

We recall that for each family $\SsChemins{n}$, $\mathrm{diam}(\SsChemins{n})= \max d_1(\mathbf{S},\mathbf{T})
\leq n(n+1)/2$.
Hence, combining Proposition \ref{Lem:Courbure}with Eq. \eqref{Eq:Mixing} gives our main result:
\begin{theo}\label{Th:Mix}
For meanders, excursions and path with walls, Algorithm \ref{Algo:CM} returns an almost uniform sample
of $\pi$, as soon as $T \gg n^3$. Precisely, for any itinialization of Algorithm \ref{Algo:CM},
$$
\parallel \mathbf{S}(T) -\pi\parallel_{\mathrm{TV}} \leq \mathrm{diam}(\SsChemins{n})(1-\kappa)^T
\leq \frac{n(n+1)}{2}\exp\left(-\frac{3}{2n^2(n+1)}T\right).
$$
\end{theo}
Another formulation of this result is that the mixing time of the associated Markov chain, defined
as usual by
\begin{equation}\label{Eq:tmix}
t_{\mbox{mix}}:=\set{\inf\ t\geq 0\ ;\ \sup_{v\in V} \parallel P^t\delta_v -\pi\parallel_{\mathrm{TV}}\leq e^{-1}}
\end{equation}
($e^{-1}$ is here by convention), is smaller than $n^2(n+1)\log n$. 
For culminating paths, the argument of Case 1c fails and \eqref{Eq:Accroissement} does not hold, we are not 
able to prove such a result as Theorem \ref{Th:Mix}. However, it seems empirically that the mixing time is also of order 
$n^3\log n$ (with a constant strongly dependent on $a,b$). A way to prove this could be the following observation: take 
$(\mathbf{S}^0,\mathbf{T}^0)=(\mathbf{S},\mathbf{T})$
two any culminating paths, and define 
$$
(\mathbf{S}^{t+1},\mathbf{T}^{t+1})=(\phi(\mathbf{S}^t,I_t,\eps_t,\delta_t),\phi(\mathbf{T}^t,I_t,\eps_t,\delta_t)),
$$
where $I_t,\eps_t,\delta_t$ are those in Algorithm \ref{Algo:CM}. The sequence 
$
\left(\parallel \mathbf{S}^t-\mathbf{T}^t\parallel_\infty \right)_t
$
is decreasing throughout the process. Unfortunately we cannot get a satisfactory bound for the time needed 
for this quantity to decrease by one.



\subsection{Related works}\label{Sec:Related}

Bounding mixing times via a contraction property over the transportation metric is quite a standard technique, the main ideas dating back to Dobrushin (1950's). A modern introduction is made in \cite{Mix}. 
For geodesic spaces, this technique has been developped in \cite{Bub} under the name
\emph{path coupling}.

As mentioned in the introduction, the Markov chain $P$ on lattice paths with uniform weights
$p_i=1/n$ has in fact already been
introduced for paths starting and ending at zero (sometimes called \emph{bridges})
in \cite{LRS}, and its mixing time has been estimated in \cite{Wilson}. 
Wilson also proves a mixing time of order $n^3\log n$, by showing that \eqref{Eq:Courbure}
holds with a different distance (namely, a kind of Fourier transform of the heights of the 
paths)\footnote{Notice that $a,b$ do not have the same meaning in Wilson's paper:
$a$ (resp. $b$) stands for the number of positive (resp. negative) steps.}. 
This is the concavity of this Fourier transform which gives a good mixing time, exactly as the concavity
of our $p_i$'s speeds up the convergence of our chain.

Wilson's method is developped only
for bridges in \cite{Wilson} and it is not completely straightforward to use it when the endpoints
are not fixed. For instance, take $n=7$ and $a=b=1$, and consider the paths $+++--++$ and $---++--$. There are more
"bad moves" (moves that take away these paths) than "good moves".

\section{\emph{Coupling From The Past} with $P$}\label{Sec:ProppWilson}

Propp-Wilson's Coupling From The Past (CFTP) \cite{PW} is a very general procedure for the exact sampling of the  stationary distribution of a Markov chain. It is efficient if the chain is monotonous with respect to
a certain order relation $\preceq$ on the set $V$ of vertices, with two extremal points 
denoted $\hat{0},\hat{1}$ (\ie such that $\hat{0}\preceq x\preceq\hat{1}$ for any vertex $x$). This is the case
here for each family $\Culmis{n}$,$\Walls{n}$,$\Exc{n}$,$\Meandres{n}$ , with the partial order
$$
\mathbf{S}\preceq \mathbf{T} \mbox{ iff } S_i\leq T_i \mbox{ for any }i.
$$
For the family $\Meandres{10}$ with $a=1,b=-2$ for instance, we have
\begin{align*}
\hat{0}=\hat{0}_{\mbox{meanders}}&=(1,1,-2,1,1,-2,1,1,-2,1),\\
\hat{1}=\hat{1}_{\mbox{meanders}}&=(1,1,1,1,1,1,1,1,1,1).
\end{align*}

We describe CFTP, with our notations, in Algorithm \ref{Algo:CFTP}. 

\begin{algorithm}
\caption{CFTP: Exact sampling of a path in $\SsChemins{n}$}
\label{Algo:CFTP}
\begin{algorithmic}
\STATE $\mathbf{S}\leftarrow\hat{0}$, $\mathbf{T}\leftarrow\hat{1}$
\STATE $\dots,I_{-2},I_{-1}\leftarrow$ i.i.d. r.v. with law $\mathbf{p}$
\STATE $\dots,\eps_{-2},\eps_{-1}\leftarrow$ i.i.d. uniform r.v. in $\set{\uparrow,\downarrow}$
\STATE $\dots,\delta_{-2},\delta_{-1}\leftarrow$ i.i.d. uniform r.v. in $\set{+,-}$
\STATE $\tau=1$
\REPEAT
  \STATE $\mathbf{S}\leftarrow\hat{0}$, $\mathbf{T}\leftarrow\hat{1}$
  \FOR{$t=-\tau$ to $0$}
  \STATE {\bf if }$\Flip{S}{I_t}{\eps_t}$ is in $\SsChemins{n}$ {\bf then } $\mathbf{S}\leftarrow  \Flip{S}{I_t}{\eps_t}{\delta_t}$
  \STATE {\bf if }$\Flip{T}{I_t}{\eps_t}$ is in $\SsChemins{n}$ {\bf then } $\mathbf{T}\leftarrow  \Flip{T}{I_t}{\eps_t}{\delta_t}$
  \ENDFOR
  \STATE $\tau\leftarrow 2\tau$
\UNTIL{$\mathbf{S}=\mathbf{T}$}
\end{algorithmic}
\end{algorithm}

We refer to (\cite{Hagg},Chap.10) for a very clear introduction to CFTP, and we only outline here the reasons why this
indeed gives an exact sampling of the stationary distribution.
\begin{itemize}
\item The output of the algorithm (if it ever ends!) is the state of the chain $P$ that has been running "since time $-\infty$",
and thus has reached stationnarity.
\item The exit condition $\mathbf{S}=\mathbf{T}$ ensures that it is not worth running the chain from $T$ steps earlier,
since the trajectory of any lattice path $\hat{0}\preceq \mathbf{R}\preceq\hat{1}$ is "sandwiched" between those of 
$\hat{0},\hat{1}$, and therefore ends at the same value.
\end{itemize}

\begin{figure}[h]
\begin{center}
\includegraphics[width=65mm]{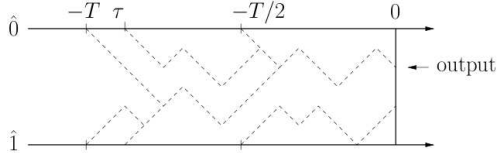}
\caption{A sketchy representation of CFTP : trajectories starting from $\hat{0},\hat{1}$ at time
$-T/2$ don't meet before time zero, while those starting at time $-T$ do.}
\label{Fig:ProppWilson}
\end{center}
\end{figure}

\begin{prop}\label{Prop:CFTP}
Algorithm \ref{Algo:CFTP} ends with probability $1$ and returns an exact sample of the uniform 
distribution over $\SsChemins{n}$.  For the families $\Walls{n}$,$\Exc{n}$,$\Meandres{n}$, 
this takes on average $\mathcal{O}(n^3(\log n)^2)$ time units.
\end{prop}
Let us mention that in the case where the mixing time is not rigorously known, Algorithm \ref{Algo:CFTP} (when it ends)
outputs an exact uniform sample and therefore is of main practical interest compared to MCMC.

\begin{proof}[Proof of Proposition \ref{Prop:CFTP}]
It is shown in \cite{PW} that Algorithm \ref{Algo:CFTP} returns an exact sampling 
in $\mathcal{O}(t_{\mbox{mix}}\log H)$ runs of the chain, where $t_{\mbox{mix}}$ is defined in \eqref{Eq:tmix}
and $H$ is the length of the
longest chain of states between $\hat{0}$ and $\hat{1}$. It is a consequence of the proof of Lemma \ref{Lem:Geodesique} 
that $H=\mathcal{O}(n^2)$. We have seen that $t_{\mbox{mix}}=\mathcal{O}(n^3\log n)$. (Recall that each test in 
Algorithm \ref{Algo:CFTP} takes, on average, $\mathcal{O}(1)$ time units.)
\end{proof}
We recall that CFTP has a major drawback compared to MCMC. For the algorithm to be correct, we have to reuse
the same random variables $I_t,\eps_t,\delta_t$, so that space-complexity is in fact linear in $n^3(\log n)^2$. 
This may become
an issue when $n$ is large. 

\section{Concluding remarks and simulations}

\begin{figure}[h!]
\begin{center}
\includegraphics[width=140mm]{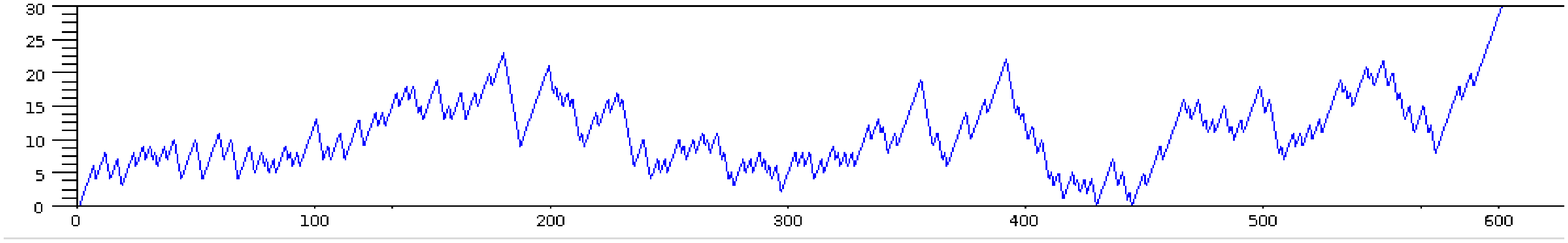}\\
\includegraphics[width=140mm]{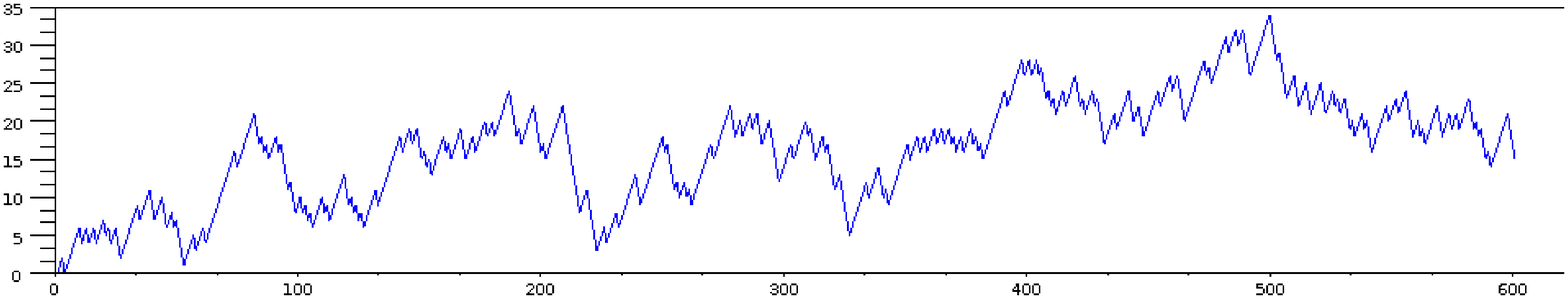}\\
\includegraphics[width=140mm]{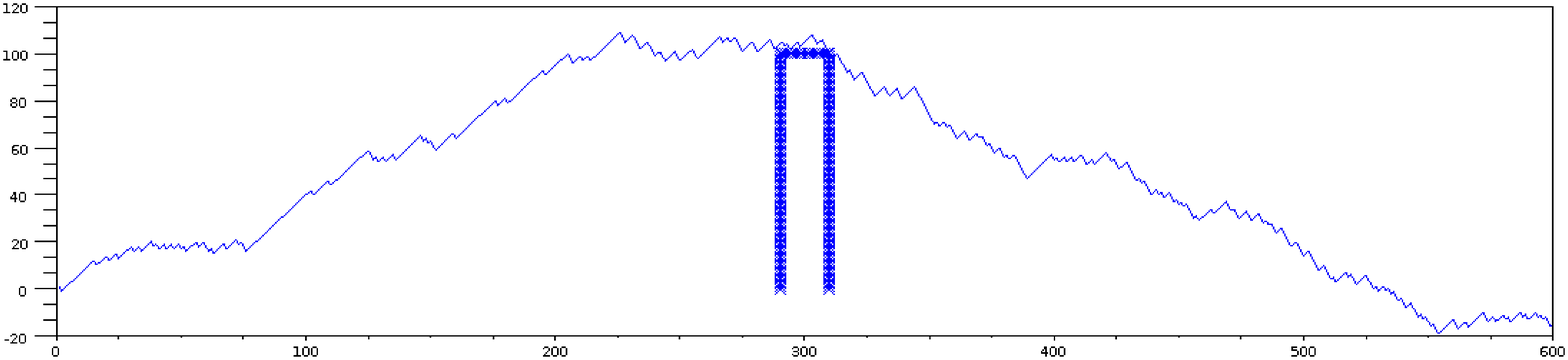}\\
\caption{(Almost) uniform paths of length $600$, with $a=1,b=2$. From top to bottom:
a culminating path, a meander, a path with wall (shown by an arch).}
\label{Fig:Simus}
\end{center}
\end{figure}

\noindent{\bf 1.} In Fig.\ref{Fig:Simus}, we show simulations of the three kinds of paths, 
for $a=1,b=2,n=600$. 
We observe that the final height of the culminating path is very low (about $30$), it would 
be interesting to use our algorithm to investigate the behaviour of this height when $n\to\infty$ ; this question was left open
in \cite{MBM}.
\vspace{4mm}

\noindent{\bf 2.} One may wonder to what extent this work applies to other families $\SsChemins{n}$
of paths. The main assumption is that the family of paths should be a geodesic space w.r.t. distance $d_1$.
This is true for example if the following condition on $\SsChemins{n}$ is fulfilled:
$$
\left(R,T\in\SsChemins{n} \mbox{ and }R\preceq S\preceq T\right) \Rightarrow S\in\SsChemins{n}.
$$
Notice however that this is quite a strong requirement, and it is not verified for culminating paths for instance.
\vspace{4mm}

\noindent{\bf 3.} A motivation to sample random paths is to
make and test guesses for some functionals of these paths, taken on average over $\SsChemins{n}$.
Consider a function $f:\SsChemins{n}\to\mathbb{R}$, we want an approximate value of
$\pi(f):=\card(\SsChemins{n})^{-1}\sum_{s\in\SsChemins{n}}f(s)$, if the exact value is out of reach by calculation. We estimate this quantity by
\begin{equation}\label{Eq:Chapeau}
\hat{\pi}(f):=\frac{1}{T} \sum_{t=1}^T f \left(\mathbf{S}(t)\right),
\end{equation}
(recall that $S(t)$ is the value of the chain at time $t$).
For Algorithm \ref{Algo:CM} to be efficient in practice, we have to bound
\begin{equation}\label{Eq:AMajorer}
\mathbb{P}\left(\left|\pi(f)-\hat{\pi}(f)\right|>r\right),
\end{equation}
for any fixed $r>0$, by a non-asymptotic (in $T$) quantity. This can be done with (\cite{JouOlli}, Th.4-5),
in which one can find concentration inequalities for \eqref{Eq:AMajorer}. The sharpness of these 
inequalities depends on $\kappa$ and on the geometrical structure of $\SsChemins{n}$.
\vspace{4mm}

\noindent{\bf Aknowledgements.} Many thanks to Fr\'ed\'erique Bassino and the other members of \textsc{Anr Gamma} for the support ; I also would like to thank \'Elie Ruderman for the English corrections.
A referee raised a serious error in the first version of this paper, I am grateful to them.

\newpage

\end{document}